\documentclass{amsart}
\usepackage{graphicx}

\usepackage{comment}
\usepackage{amssymb}
\usepackage{amsthm}
\usepackage{amsmath}
\usepackage{tikz}
\usepackage[normalem]{ulem}

\usepackage{pgfplots}
\pgfplotsset{compat=1.18}
\usepackage{tikz-3dplot}

\newtheorem{thm}{Theorem}[section]

\newtheorem{rmk}[thm]{Remark}

\newtheorem{exc}[thm]{Exercise}
\newtheorem{prob}[thm]{Problem}

\newcommand{\todomarginpar}[1]{}

\title{What Grassmann Knew: Incidence Theorems on Cubics}
\author{Will Traves}
\address{Department of Mathematics, United States Naval Academy}
\email{traves@usna.edu}

\thanks{
\noindent\textbf{Keywords}: synthetic geometry, cubic curve, incidence, line arrangement, Grassmann-Cayley algebra}  \thanks{\noindent\textbf{2020 Mathematics Subject Classification}: Primary: 
14H50, 15A75, 51A20;  Secondary: 14N20, 52C35}

\begin{document}

\maketitle

\begin{abstract}
Traves and Wehlau \cite{TW} recently gave a straightedge construction that checks whether 10 points lie on a plane cubic curve. They also highlighted several open problems in the synthetic geometry of cubics. Hermann Grassmann investigated incidence relations among points on cubic curves in three papers \cite{G2,G3,G1} appearing in Crelle's Journal from 1846 to 1856. Grassmann's methods give an alternative way to check whether 10 points lie on a cubic. Using Grassmann's techniques, we solve the synthetic geometry problems introduced by Traves and Wehlau. In particular, we give straightedge constructions that find the intersection of a line with a cubic, find the tangent line to a cubic at a given point, and find the third point of intersection of this tangent line with the cubic. As well, given 5 points on a conic and a cubic and 4 additional points on the cubic, a straightedge construction is given that finds the sixth intersection point of the conic and the cubic. The paper ends with two open problems.
\end{abstract}

\section{Introduction}

In our paper, Ten Points on a Cubic \cite{TW}, David Whelau and I gave a straightedge construction to check whether 10 points in the plane lie on a cubic curve. We ended our paper with three open problems in synthetic geometry. In the meantime, I learned that the great German mathematician Hermann Grassmann did extensive work on this problem in between 1846 and 1856, arriving at a completely different construction that checks whether 10 points lie on a cubic. Grassmann's work was popularized by Alfred North Whitehead in his book Universal Algebra and championed by the American physicist David Hestenes. In the next section we describe Grassmann's approach to projective geometry via incidence structures and illustrate his methods using conics. In Section \ref{section:cubics} we give Grassmann's description of cubic curves. In the following section we use Grassmann's ideas to solve the three synthetic geometry problems posed in Traves and Wehlau \cite{TW}. We go further, giving a straightedge construction that checks whether a point on a cubic defined by nine points is a flex. We close with some pointers to the literature and a few more open problems. 

\section{Projective Geometry from Incidence Structures} \label{section:incidence}

Hermann Grassmann developed a remarkable algebraic system that models a vast array of geometric settings. One descendent of this system is the Grassmann-Cayley algebra, which describes geometric operations with linear subspaces. This algebra is usually presented in terms of the exterior algebra but it admits a simplified version in the setting of the projective plane that encodes straightedge constructions. 

Points in $\mathbb{P}^2$ are represented by equivalence classes of nonzero vectors with three components: the point $[x_0:x_1:x_2]$ is equivalent to each point of the form $[kx_0: kx_1: kx_2]$, where $k$ is a nonzero scalar. The point $(x_1,x_2)$ in the Euclidean plane is usually identified with the class $[1:x_1:x_2]$. Similarly, lines in $\mathbb{P}^2$ are represented by a point in a dual $\mathbb{P}^2$: the line $L_0x_0 + L_1x_1 + L_2x_2 = 0$ is represented by the point $[L_0:L_1:L_2]$. It will be convenient to introduce a zero-point $[0:0:0]$ and a zero line $0x_0 + 0x_1 + 0x_2 = 0$ as degenerate objects for book-keeping purposes. 

Our two operations are the meet and the join. Given two lines $L$ and $M$ in $\mathbb{P}^2$ with equations $L_0x_0 + L_1x_1 + L_2x_2 = 0$ and $M_0x_0 + M_1x_1 + M_2x_2 = 0$, respectively, the meet $LM$ of the line $L$ with the line $M$ is the point  $ [N_0: N_1: N_2]$ given by the cross product of the coefficient vectors of $L$ and $M$: 
$$  \begin{array}{lll} [N_0: N_1: N_2] & = & [L_0: L_1: L_2] \times [ M_0: M_1: M_2] \\
& = & [L_1M_2-L_2M_1: L_2M_0-L_0M_2 : L_0M_1 - L_1M_0].\end{array}$$
The reader should check that the meet operation is well-defined, though if $L=M$ we obtain the zero point. 

The join of two points is defined by dualizing the meet of two lines. If $p = [p_0:p_1:p_2]$ and $q = [q_0:q_1:q_2]$ are two points, their join is the line $pq$ through the two points, whose equation is $L_0x_0 + L_1x_1 + L_2x_2 = 0$ where $$  \begin{array}{lll} [L_0: L_1: L_2] & = & [p_0: p_1: p_2] \times [ q_0: q_1: q_2] \\
& = & [p_1q_2-p_2q_1: p_2q_0-p_0q_2 : p_0q_1 - p_1q_0].\end{array}$$
The join of a point with itself is the zero line. 

We supplement the cross product with the dot product to define an operator, the bracket, whose vanishing encodes incidence relations. The bracket $(p.q.r)$ of three points $p$, $q$, and $r$ is the scalar triple product: $$(p.q.r) = p \cdot (q \times r).$$ Though this is not well-defined as a number, whether the value is zero or not is well-defined and this is sufficient for our purposes. It will be convenient to introduce an equivalence relation on our scalars: all nonzero scalars are in one class and the zero scalar is in another class. For two scalars $a$ and $b$, we write $a \equiv b$ to indicate that $a$ and $b$ are in the same class. 
Note that $(p.q.r)=0$ means that either some of the three points coincide or the three points are distinct and collinear. We say that the points $p$, $q$ and $r$ are the three factors of the bracket $(p.q.r)$. Similarly, we define the bracket $(L.M.R)$ of three lines $L$, $M$, and $R$ as $L \cdot (M \times R)$. This bracket is zero when two of the lines are identical or when all three lines meet in a common point. 

If $L=pq$ is a line and $r$ is a point, the expression $Lr$ is defined to be a scalar that vanishes when $(p.q.r)=0$, that is, $Lr=0$ when $r$ lies on $L$. Similarly, the point $p$ lies on the line $M = qr$ if and only if $pM := (p.q.r) = 0$. If all three points are distinct, the scalar given by the expression $pqr$ is equivalent to any of its permutations but this is not true for more general expressions. For instance, consider the expression 
$pqrs$. This really ought to have parentheses to indicate the order of operations since $(pq)(rs)$ is the point of intersection of the lines $pq$ and $rs$, while $((pq)r)s = (p.q.r)s$ is a multiple of the point $s$ (i.e. either $s$ or the zero point). Since these are different in general we will evaluate all expressions from left to right unless indicated otherwise by parentheses. We also follow Grassmann's practice and use periods to indicate different geometric items, so $pq.rs$ represents the point $(pq)(rs)$. 

As practice with the notation, we describe the expression $abAcBd$. Throughout the paper lower case letters denote points and upper case letters denote lines. So this expression requires us to first intersect the line $ab$ with the line A, then join the resulting point with the point $c$ to obtain a line that we intersect with the line $B$ and then join the resulting point with the point $d$ to obtain a line. The resulting line could be the zero line. This happens when $a=b$, $ab=A$, $abA=c$, $abAc=B$ or $abAcB=d$.   

When $x = [x_0:x_1:x_2]$ is a variable point in $\mathbb{P}^2$, and $p$ and $q$ are two distinct points in $\mathbb{P}^2$, the expression $(x.p.q)=0$ is the equation of the line $pq$. The equation $xaAbBcx=0$ is another familiar object: the equation of a conic. The expression contains two copies of $x$ so expanding using the cross and scalar product gives a degree-2 equation in the entries of $x$. 

\begin{exc} (a) Show that if the points $a$, $b$, and $c$ are not collinear, none of them lie on $A$ or $B$ and $A$ and $B$ are distinct lines, then the conic with equation $xaAbBcx=0$ passes through the 5 points $a$, $b$, $AB$, $abB$, and $bcA$. \\
(b) Find choices of the parameters so that $xaAbBcx=0$ is the equation of two crossing lines. \\
(c) Find choices of the parameters so that $xaAbBcx=0$ is the square of a linear form.\\
(d) Pascal's Theorem says that if points $a$, $b$, $c$, $a_1$, $b_1$, and $c_1$ lie on a conic then the points $(ab_1)(a_1b)$, $(ac_1)(a_1c)$ and $(bc_1)(b_1c)$ are collinear. Use the relation $(ab)(cd) = (a.b.d)c-(a.b.c)d$ to deduce Pascal's Theorem from Grassmann's description of a conic.  [Note that the two points in this relation are well-defined, even if we scale the points $a$, $b$, $c$, and $d$.] 
\end{exc}

\section{Grassmann on Cubics} \label{section:cubics}

Grassmann proposed to parameterize a cubic curve $Q$ using 6 points, $a$, $a_1$, $b$, $b_1$, $c$ and $k$, and 3 concurrent lines $A$, $B$, and $C$: 
$$Q = \{ x \in \mathbb{P}^2: \, (xaAa_1.xbBkCb_1.xc)=0 \}.$$
We preserve Grassmann's idiosyncratic nomenclature.  The point $x$ lies on the cubic $Q$ if the three lines $L=xaAa_1$, $M=xbBkCb_1$, and $R=xc$ meet in at least one common point or if any of the constructions required to produce $L$, $M$ and $R$ are degenerate. Figure \ref{fig:cubicconstruction} on the next page illustrates the construction. 

\newcommand{\intregion}[1]{
\begin{tikzpicture}[scale=#1]
 
\filldraw[black] (3.55,2.01) circle (4pt) node[anchor=south east]{$x$};
\filldraw[black] (1.04,0.65) circle (4pt) node[anchor=north west]{$a$};
\filldraw[black] (6.05,3.37) circle (4pt) node[anchor=south east]{$b$};
\filldraw[black] (9.69,1.67) circle (4pt) node[anchor=west]{$k$};
\filldraw[black] (5.32,-1.54) circle (4pt) node[anchor=north west]{$b_1$};
\filldraw[black] (2.96,-0.41) circle (4pt) node[anchor=south east]{$c$};
\filldraw[black] (1,-2) circle (4pt) node[anchor=north]{$a_1$};

\node at (-2.3,2.5) {$A$};
\node at (10.5,5.8) {$B$};
\node at (11.5,0.15) {$C$};

\draw[ultra thick,black] (-1.24,-4.02) -- (-3,3);
\draw[ultra thick,black] (2.29,-3.12) -- (3.95,3.64);
\draw[ultra thick,black] (1.6,-2.76) -- (10,0);
\draw[ultra thick,black] (9.08,5.02) -- (10,0);
\draw[ultra thick,black] (7.31,3.26) -- (11.1,7.02);
\draw[ultra thick,black] (-2,-1) -- (9.08,5.02);
\draw[ultra thick,black] (4.05,-3.02) -- (-2,-1);
\draw[ultra thick,black] (12.55,1.16) -- (7.48,-1.15);
\end{tikzpicture}}

\begin{figure}[h!t]
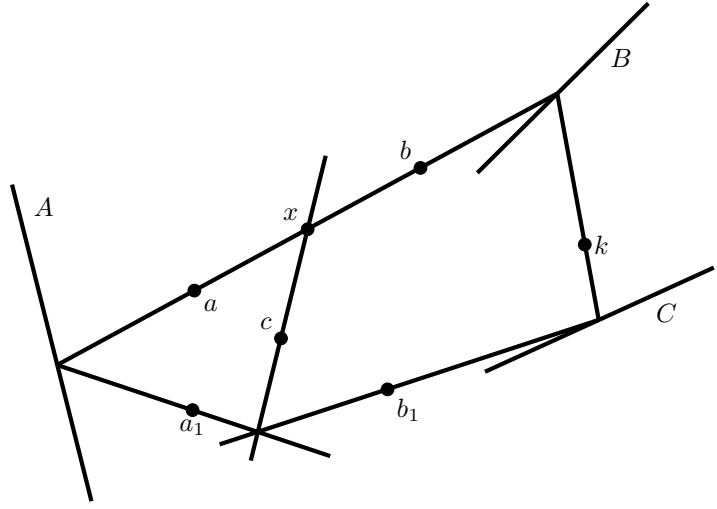

\begin{center}
\intregion{0.6}
\end{center}
\caption{The point $x$ lies on the cubic $(xaAa_1.xbBkCb_1.xc)=0$.}
\label{fig:cubicconstruction}
\end{figure}

The equation $(xaAa_1.xbBkCb_1.xc)=0$ defining the curve $Q$ contains the point $x = [x_0:x_1:x_2]$ three times, so $Q$ is the zero set of a degree-3  polynomial in the variables $x_0$, $x_1$, and $x_2$. The general cubic equation is 
$$ \resizebox{0.85\width}{0.85\height}{$ a_{0}x_0^3+a_{1}x_0^2x_1+a_{2}x_0^2x_2+a_{3}x_0x_1^2+a_{4}x_0x_1x_2+a_{5}x_0x_2^2+a_{6}x_1^3+a_{7}x_1^2x_2+a_{8}x_1x_2^2+a_{9}x_2^3 = 0$},$$
so each cubic corresponds to a point $[a_0:a_1:\cdots:a_9]$ in $\mathbb{P}^9$. Requiring the cubic to pass through a point $p=[p_0:p_1:p_2]$ imposes a linear constraint on the coefficients of the cubic. Requiring the cubic to pass through 9 points in general position uniquely determines the cubic. 

Given nine points $a$, $b$, $c$, $d$, $e$, $f$, $g$, $h$ and $i$ in general position, we show how to choose Grassmann's parameters to produce the cubic through the given points. We set $A = de$, $B = ef$, and $a_1 = af.cd$. Then we define six new points $$ \begin{array}{ll}
g_1 =gaAa_1.gc, & g_2 = gbB \\
h_1 =haAa_1.hc, & h_2 = hbB \\
i_1 =iaAa_1.ic, & i_2 = ibB. \end{array}$$
The line $C$ is defined as $ei_1$ and we define two additional points $y = h_1g_1Cg_2.fh_1$ and $z=g_1h_1Ch_2.fg_1$, the line $K = yz$ and the points $k= K.i_1i_2$ and $b_1 = kg_2Cg_1.kf$. Then the cubic $Q$ defined by $(xaAa_1.xbBkCb_1.xc)=0$ turns out to pass through the nine given points $a$, $b$, $\ldots$, $i$. 

Grassmann's original argument \cite{G3} justifying this construction is subtle; it involves showing that an associated cubic is the union of three lines. Alfred North Whitehead translated Grassmann's argument into English for his book, Universal Algebra \cite{Whitehead}. Whitehead added a few details and tried to streamline Grassmann's proof, but his version contains a logical flaw. Whitehead assumes that every cubic can be represented by an appropriate choice of the parameters and comes up with a construction of the parameters under that assumption. Strangely, Hestenes and Ziegler \cite{HZ} repeat Whitehead's proof for this construction almost word for word. We present a direct proof that is very close to Grassmann's original proof. 

\begin{thm} \label{thm:ninepts}
When the nine points $a$, $b$, $c$, $d$, $e$, $f$, $g$, $h$ and $i$ lie in linearly general position (no three on a line), the construction described above produces a cubic $(xaAa_1.xbBkCb_1.xc)=0$ that passes through the nine given points.   
\end{thm}

\begin{proof}
Setting $L = xaAa_1$, $M = xbBkCb_1$, and $R = xc$, the equation for the cubic becomes $(L.M.R)=0$. We see that $L=0$ when $x=a$, $M=0$ when $x=b$, and $R=0$ when $x=c$ so the points $a$, $b$, and $c$ lie on the cubic. When $x=d$, $daA = d$ so $L \equiv da_1 \equiv cd \equiv R$ so $(L.M.R)=0$. When $x=e$, $L \equiv ea_1$, $M \equiv eb_1$, and $R \equiv ec$ share the common point $e$ so $e$ lies on the cubic. When $x=f$, $L\equiv fa$, $M \equiv fkCb_1 \equiv fk$ (because $b_1$ lies on $kf$), and $R\equiv fc$ share the common point $f$ so $f$ lies on the cubic. When $x=g$, $L.R = g_1$ and since $b_1$ lies on $g_2kCg_1$ the three points $b_1$, $g_2kC$ and $g_1$ are collinear and $g_1$ lies on $M = g_2kCb_1$. So when $x=g$ the three lines $L$, $M$ and $R$ meet in $g_1$ so $g$ lies on the cubic. Similiarly, when $x=i$, the three lines $L$, $M$ and $R$ meet in $i_1$ so $i$ lies on the cubic.  

It remains to consider the case $x=h$. We first establish a preliminary result: the point $b_1 = g_2kCg_1.kf$ lies on the line $kh_2Ch_1$. The proof involves showing that an auxiliary cubic $(xf.xg_2Cg_1.xh_2Ch_1)=0$ contains $k$. So we turn our attention to this auxiliary cubic. It turns out that this auxiliary cubic is the union of three lines. If $x$ lies on $B=ef$ then $xg_2 = xh_2 = xf = B$ so each of the lines $xf$,  $xg_2Cg_1$, and $xh_2Ch_1$ contain $BC = e$ and $(xf.xg_2Cg_1.xh_2Ch_1)=0$; so the auxiliary cubic contains the line $B$. If $x$ lies on $C$ then each of the lines $xf$, $xg_2Cg_1$ and $xh_2Ch_1$ contain $x$, so the auxiliary cubic also contains the line $C$. We show that the remaining component of the auxiliary cubic is the line $K$ containing $y$ and $z$.   

First we show that both $y$ and $z$ lie on the cubic. Since $y=h_1g_1Cg_2.fh_1$, $(yfh_1)=0$ and $h_1$ lies on $yf$. The vanishing of the expression $h_1g_1Cg_2y$ means that $h_1$ lies on the line $yg_2Cg_1$. Clearly, $h_1$ lies on $yh_2Ch_1$, so $h_1$ lies on each of the factors of the auxiliary cubic evaluated at $x=y$, showing that $(yf.yg_2Cg_1.yh_2Ch_1)=0$. This guarantees that $y$ lies on the auxiliary cubic. Similarly, one can show that $g_1$ lies on each of the factors of the auxiliary cubic evaluated at $x=z$ so $z$ also lies on the auxiliary cubic. 

Now we argue that neither $y$ nor $z$ lie on the lines $B$ and $C$. To see that $y$ cannot lie on $B=ef$, note that $h_1$ lies on $yf$ so if $y \in B$ then $h_1 \in B$. But then $c \in h_1y = B$ and so $(cef)=0$. This violates the assumption that the nine given points are in linearly general position, so we conclude that $y \not\in  B$. If $y \in C$ then $g_2 \in C$. But $g_2 \in B$ by construction, so $g_2 = BC = e$ lies on $gb$, violating the assumption of linearly general position. So $y \not\in C$. Similar arguments show that $z$ does not lie on $B$ or $C$. 

Since neither $y$ nor $z$ lie on the two components $B$ and $C$ of the auxiliary cubic, the points $y$ and $z$ must both lie on the remaining component, the line $K=yz$. It follows that any point on the line $K$ also lies on the auxiliary cubic. In particular, the point $k = i_1i_2.K$ lies on the auxiliary cubic:
$$ (kf.kg_2Cg_1.kh_2Ch_1)=0.$$
The three lines $kf$, $kg_2Cg_1$ and $kh_2Ch_1$ meet in the point $b_1 = kg_2Cg_1.kf$. 

Now we check that $x=h$ lies on the original cubic $(xaAa_1.xbBkCb_1.xc)=0$. Note that the middle factor evaluated at $x=h$ is a line through $h_1$ since 
$$ hbBkCb_1h_1 = h_2kCb_1h_1 = kh_2Ch_1b_1 = 0, $$  
as $b_1$ lies on the line $kh_2Ch_1b_1$. The point $h_1 = haAa_1.hc$ also lies on the first and last factors in the cubic when evaluated at $x=h$. The three factors of the cubic all pass through $h_1$ when $x=h$, so $h$ lies on the cubic. 
\end{proof} 

\section{Synthetic Geometry}

In this section, we solve some synthetic geometry problems using Grassmann's description of cubic from Section \ref{section:cubics}. Several of these problems were first introduced in Traves and Wehlau \cite{TW}. 

\begin{thm} \label{thm:intlinecubic} Given six points and three lines in general linear position defining the cubic $(xaAa_1.xbBkCb_1.xc)=0$, the third point of intersection of the line through $a$ and $b$ with the cubic is $(abAa_1.abBkCb_1)c.ab$.  
\end{thm}

\begin{proof} 
As motivation,   we first show that the formula for the third point follows from basic geometric arguments. Let $y = \mu a + \nu b$ lie on $ab$ and suppose that it lies on the cubic. Then $(yaAa_1.ybBkCb_1.yc)=0$. Note that $ya = (\mu a + \nu b)a = \nu ba$ and $yb = (\mu a + \nu b)b = \mu ba$. Substituting into the first and second factors, we find $\mu\nu(abAa_1.abBkCb_1.yc)=0$. Setting $p = abAa_1.abBkCb_1$ we get $\mu\nu(p.yc)=0$. So either $\mu=0$ (when $y=a$), $\nu=0$ (when $y=b$), or $(pyc)=0$, when $y$ lies on both $pc$ and $ab$, so $y=pc.ab = (abAa_1.abBkCb_1)c.ab$. 

Now we consider the sufficiency statement. The computation for $p$ cannot be degenerate, $p \equiv 0$, otherwise the points and lines defining the cubic are not in general linear position. If $pc$ passes through $a$ then $ab$ is tangent to the cubic at $a$ (see Theorem \ref{thm:tangency} below). If $pc$ passes through $b$ we can swap the labels on $a$ and $b$ to see that $ab$ must be tangent at $b$. So we assume that $y = \mu a+ \nu b$ lies on $pc$ and $\mu \nu \neq 0$. Then we can write $b = (y-\mu a)/\nu$ and $a = (y-\nu b)/\mu$. The condition $y.pc = 0$ forces $p.yc$ = 0 so 
$$ \begin{array}{lll} \frac{1}{\mu\nu} (yaAa_1.ybBkCb_1.yc) & = &  \left(\left(\frac{y-\mu a}{\nu}\right)aAa_1.\left(\frac{y-\nu b}{\mu}\right)bBkCb_1.yc\right) \\ 
& = & (baAa_1.abBkCb_1.yc) \\ & = & p.yc = 0,
\end{array}$$
and $y=pc.ab$ lies on the cubic curve. Thus $y$ is the third point of intersection of $ab$ with the cubic.
\end{proof}

\begin{rmk} A smooth cubic $Q$ with a fixed flex point, $o$, admits a group law: given points $a$ and $b$, 
$$a+b = [o,[a,b]],$$ where $[x,y]$ denotes the third point of intersection of the line $xy$ with the cubic. It follows from Theorem \ref{thm:intlinecubic} that the group operation $a+b$ on the cubic can be computed with a straightedge if all the relevant points -- $a$, $b$ and $o$ -- are given. \end{rmk}

\begin{thm} \label{thm:tangency} Given six points $a$, $b$, $c$, $a_1$, $b_1$, and $k$, and three lines $A$, $B$, and $C$ defining the cubic $(xaAa_1.xbBkCb_1.xc)=0$, the tangent to the cubic at the point $a$ is $(abBkCb_1.ac)a_1Aa$. 
\end{thm}

\begin{proof}

Our approach to finding the tangent line involves an infinitesimal argument. The tangent line at $a$ points towards another point $v$ that we locate by substituting $x=a+\epsilon v$ into the equation for the cubic and expanding terms to get 
$$ \begin{array}{ll} & \left( (a+\epsilon v)aAa_1.(a+\epsilon v)bBkCb_1.(a+\epsilon v)c\right) \\ = & \epsilon(vaAa_1.abBkCb_1.ac) + \\ & \epsilon^2(vaAa_1.abBkCb_1.vc + vaAa_1.vbBkCb_1.ac) + \\ & \epsilon^3(vaAa_1.vbBkCb_1.vc).\end{array}$$
The line $av$ is a tangent line to the cubic at $a$ when the $\epsilon$ term vanishes. That is, when the line $vaAa_1$ contains the point $p = (abBkCb_1.ac)$. This occurs when $vaAa_1p = 0$. This says that the tangent line hits $A$ at the same point as $pa_1A$. So the tangent line to the cubic at $a$ is $pa_1Aa = (abBkCb_1.ac)a_1Aa$. 
\end{proof}

\begin{rmk}
The algorithm in Theorem \ref{thm:tangency} produces a unique tangent line when the cubic is smooth at $a$ but when the cubic is singular at $a$ the tangent line produced depends on the choice of the parameters of the cubic. Computer exploration has shown that the algorithm can be made to produce any line through $a$ by varying the parameters that produce the cubic curve. This behavior should be expected since the tangent space to the cubic at the singular point is the entire plane.  
\end{rmk}

\begin{thm}
\label{thm:inttangcubic}
The third point of intersection of the tangent line to the cubic $(xaAa_1.xbBkCb_1.xc)=0$ at $a$ hits the cubic in a third point that can be determined as follows. Setting $q = (abBkCb_1.ac).a_1A$, the tangent line to the cubic at $a$ is $aq$. The conic $xbBkCb_1x=0$ passes through $b$, $b_1$, $BC$, $b_1kB$, and $bkC$; set $y$ equal to the second point where the line $cb_1CkBb$ through $b$ meets this conic. Further, set $z$ equal to the point at the intersection of the lines $b_1cCkBb$ and $b_1c$. Then the conic $(qa_1.xc.xbBkCb_1)=0$ passes through $a$, $b$, $c$, $y$ and $z$ and it meets the tangent line $aq$ in another point $w = aq.y(bz.(ab.yc)(aq.zc))$. The point $w$ is the third point of intersection of the tangent line to the cubic $(xaAa_1.xbBkCb_1.xc)=0$ at $a$ with the cubic. 
\end{thm}

\begin{proof}
We check that $y$ lies on the conic $(qa_1.xc.xbBkCb_1)=0$; indeed, we will show that when $x=y$, the middle and right factors are equal up to a scaling. We show that the right factor $ybBkCb_1$ contains both $y$ and $c$. The point $y$ satisfies $ybBkCb_1y=0$ since $y$ lies on the conic $xbBkCb_1x=0$, so the line $ybBkCb_1$ contains $y$. We know that $y$ lies on the line $cb_1CkBb$ and reading this backwards we find that $c$ lies on the line $ybBkCb_1$ too. So the line $ybBkCb_1$ equals the line $yc$.  

We note that $z$ lies on the conic $(qa_1.xc.xbBkCb_1)=0$: $$zbBkCb_1 \equiv b_1c \equiv zc, $$ because the four points $z$, $c$, $b_1$ and $zbBkC$ are collinear.

The point $w$ is on the tangent line by construction and $w$ lies on the original cubic $(xaAa_1.xbBkCb_1.xc)=0$ because $waAa_1 \equiv aqAa_1 \equiv qa_1$ and $w$ satisfies the conic equation $(qa_1.xc.xbBkCb_1)=0$. 
\end{proof}

\begin{rmk}
Another way to find the third point of intersection of the tangent line to the cubic as in Theorem \ref{thm:inttangcubic} uses the $8 \Rightarrow 9$ Theorem. Recall that this theorem says that any cubic through 8 of the 9 points of intersection of two cubics must pass through the ninth point too. Join the point of tangency $a$ to two other points $p_1$ and $q_1$. Let $p_2$ (respectively $q_2$) be the third point of intersection of $p_1a$ (respectively $q_2a$) with the cubic. Let $p_1q_1$ (respectively $p_2q_2$) meet the cubic in the third point of intersection $r_1$ (respectively $r_2$). Finally, let $r_1r_2$ meet the cubic in the third point $r_3$. 
The union of the lines $p_1p_2$ and $q_1q_2$ (both containing $a$) with $r_1r_2$ form a cubic meeting the original cubic in three more points, $r_3$ and $a$ (counted with multiplicity two). If we denote the tangent line at $a$ by $L$, the union of the three lines $p_1q_1$, $p_2q_2$, and $L$ meets the original cubic in eight points, $p_1$, $q_1$, $p_2$, $q_2$, $r_1$, $r_2$, and $a$ (counted twice\footnote{Counting in this way recalls the joke: There are three kinds of mathematicians, those who can count ... and those who can't.}). The $8 \Rightarrow 9$ Theorem forces $r_3$ to lie on the tangent line $L$, so the third point of intersection of the tangent line $L$ with the original cubic lies at $r_3$.  
\end{rmk}

\begin{rmk}
Recall that a point $a$ on a cubic is a flex if its tangent line intersects the cubic at $a$ to order 3. This is equivalent to saying that the tangent line only hits the cubic at the point $a$, or that the point $w$ defined in Theorem \ref{thm:inttangcubic} is congruent to $a$: $w \equiv a$. 
\end{rmk}

Our next theorem says that given 5 points of intersection of a cubic and conic we can find the remaining point of intersection using a straightedge (if the points are in general linear position). 

\begin{thm} \label{thm:cc} Let nine points in general linear position be given, five of the points, $a$, $c$, $d$, $e$ and $f$, on both a conic and a cubic and additional points $b$, $g$, $h$, and $i$ on the cubic. Then the conic is $\{x \in \mathbb{P}^2: xaAa_1Bcx=0\},$ where $A=de$, $B=fe$, and $a_1=af.dc$. The cubic has equation $(xaAa_1.xbBkCb_1.xc)=0$ with parameters chosen as in Theorem \ref{thm:ninepts}. Then the points $e$ and $f$ both lie on the auxiliary cubic $(xa_1Aa. xb_1CkBb.xc)=0$; set $y$ equal to the third point of intersection of the line $ef$ with the auxiliary cubic. Then the point $z = (yc.ya_1Aa)$ lies on both the conic and the original cubic. 
\end{thm}

\begin{proof}
The auxiliary cubic $(xa_1Aa.xb_1CkBb.xc)=0$ contains the common intersection point $e$ of the lines $A$, $B$, and $C$ since all three factors pass through $e$ when $x=e$. Now we show that $f$ also lies on this auxiliary cubic. Recall that $f$ lies on the original cubic so $(faAa_1.fbBkCb_1.fc)=0$. Since $a_1 = af.cd$, the points $a$, $a_1$ and $f$ are collinear so the first factor in the original cubic is $af$. This line hits the third factor $fc$ in the point $f$. We conclude that the middle factor $fbBkCb_1$ is a line passing through $f$. Reading $fbBkCb_1f=0$ backwards gives $fb_1CkBbf=0$, which shows that if $x = f$ satisfies the equation for the auxiliary cubic $(xa_1aAa.xb_1CkBb.xc)=0$ because all three factors pass through $f$. We can use Theorem \ref{thm:intlinecubic} to produce the third point of intersection $y$ of the line $B = ef$ with the auxiliary cubic.  Setting $z = yc.ya_1Aa$, we produce a sixth point of intersection of the conic and the original cubic. Note that $zaAa_1B \equiv y$ so $z$ satisfies the equation for the conic, $zaAa_1Bcz = ycz = 0$. When $x=z$ the first and last factors of the cubic $(xaAa_1.xbBkCb_1.xc)$ pass through $y$ by our definition of $z = yc.ya_1Aa$. The second factor also passes through $y$ because $zbBkCb_1y \equiv yb_1CkBbz = 0$, with the last equality following from the fact that $y$ lies on the auxiliary cubic so the middle factor $yb_1CkBb$ must pass through the intersection point $z$ of the other two factors.   This shows that $z$ is the remaining point of intersection of the original cubic and the conic. \todomarginpar{May need to check carefully that $z$ cannot be $a$, $c$, $d$, $e$ or $f$. In these cases we should violate the general linear position hypothesis.}  
\end{proof}

\begin{rmk}
We give another way to compute the sixth point of intersection of the conic and the cubic in Theorem \ref{thm:cc} using the $8 \Rightarrow 9$ Theorem. Let $cd$ meet the cubic in the new point $r$ and $ef$ meet the cubic in the new point $s$. Let $rs$ meet the cubic in the new point $t$. Then the union of $rs$ and the conic meets the cubic in 9 points, $r$, $s$, $t$, $a$, $c$, $d$, $e$, $f$, and the point $z$ that we are trying to construct. The union of the three lines $at$, $cd$, and $ef$ contains $r$ and $s$, so the union must also contain $z$ by the $8 \Rightarrow 9$ Theorem. Since the original cubic went through 9 points in general position, the cubic is irreducible and $z$ must lie on the line $at$. So we can construct $z$ as the third point of intersection of the line $at$ with the cubic. 
\end{rmk}

\todomarginpar{Add numerical examples to provide additional confidence that the results hold. Probably later.} 

\section{Extensions and Open Problems.}

We close with some open problems, pointers to further reading, and comments about how this work connects to related topics. Whitehead's book on Universal Algebra \cite{Whitehead} contains additional straightedge characterizations of cubics. Grassmann's original papers in Crelle's Journal \cite{G2,G3,G1} remain a great source of inspiration and provide an opportunity to brush up on your German language skills too. Forder \cite{Forder} gives an English language account of Grassmann's theories with applications in geometry and engineering.  Sturmfels \cite{SIT} describes the Grassmann-Cayley algebra, the higher dimensional version of Grassmann's extension theory, and explains its connection with classical invariant theory.   Hestenes and Ziegler \cite{HZ} explain how  Grassmann's theory of extension connects to the larger theory of Geometric Algebra. Bayro-Corrochano \cite{BC} gives a detailed account of Geometric Algebra with a focus on applications in computer vision, graphics, and neurocomputing. Computing with Grassmann expressions gets easier if you use a computer algebra package, such as John Browne's Grassmann algebra package for Mathematica \cite{Browne}.

\begin{prob} \label{fourandline} Suppose that we are given 9 points $p_1,\ldots,p_9$ in the plane lying on a single cubic. Given another point $p_0$ such that $p_0$, $p_1$, $p_2$, $p_3$ and $p_4$ lie on a single conic, give a straightedge construction that produces the line through the two remaining points of intersection of the cubic and the conic, or show that no such construction exists.
\end{prob} 

\begin{rmk} It is tempting to try to use the $8 \Rightarrow 9$ Theorem to construct two points on the line through the two remaining points of intersection in Problem \ref{fourandline} but it seems that it is only possible to construct one point on the line using this method.  
\end{rmk}

\begin{rmk}
Problem \ref{fourandline} admits a solution if we allow ourselves to work in 3-space, but it is an open problem to find a solution that only uses linear operations in $\mathbb{P}^2$. Here is a sketch of the solution in 3-space. Let $s$ and $t$ denote the two unknown points of intersection of the conic and the cubic. Note that it is possible to construct points $P_i$ and $Q_i$ ($i \in \{1,2,3,4\}$) on the cubic so that $P_i$, $Q_i$, $s$, $t$, $p_1$ and $p_2$ all lie on a conic. After change of coordinates, we can assume that $p_1 = [0:0:1]$ and $p_2 = [0:1:0]$ so that the conics all have the form $a_0x_0^2 + a_1x_0x_1 + a_2x_0x_2 + a_3x_1x_2 = 0.$ That is, each conic corresponds to a point $[a_0:a_1:a_2:a_3] \in \mathbb{P}^3$. In fact, requiring the $i^\text{th}$ conic to pass through $P_i$ and $Q_i$, we find that the conic corresponds to a point on a line in $\mathbb{P}^3$. Define the injective map $F: \mathbb{P}^2 \rightarrow \mathbb{P}^3$ by $F([x_0:x_1:x_2]) = [x_0^2: x_0x_1: x_0x_2: x_1x_2])$, whose image is a quadratic surface. Then the $i^\text{th}$ conic corresponds to a point on the line joining $F(P_i)$ and $F(Q_i)$. There exists a \emph{unique} line transversal to all four of these lines and this transversal can be constructed by an appropriate combination of meets and joins. This line hits the quadratic surface at two points $F(s)$ and $F(t)$ and we can recover the points $s$ and $t$ from Von Staudt constructions. Though the details of this sketch can be filled in, it is not clear how to push the entire construction down into $\mathbb{P}^2$. 
\end{rmk}

\begin{prob} There are many possible formats for a degree-4 analogue of Grassmann's parameterization of degree-3 curves. For instance, 
$$(xa_1Aa_2Ba_3. xb_1Cb_2Dx. xc_1Ec_2Fc_3) = 0$$ parameterizes a degree-4 curve since the variable point $x$ appears four times in the expression. Another possible parameterization might be 
$$(xa_1Aa_2. xb_1Bb_2Cb_3Dx. xc_1Ec_2Fc_3) = 0.$$ A tangent-space computation shows that the first expression parameterizes almost all degree-4 curves. To show that an expression of this form can be used to describe the points on any degree-4 curve, we need an analogue of Theorem \ref{thm:ninepts}, showing that for every set of 14 points there is a choice of the 14 parameters $a_1,A,\ldots$ so that the parameterized curve passes through all 14 points. Such a result would immediately give a straightedge construction that checks whether 15 points lie on a degree-4 curve. 
\end{prob}

\section{Acknowledgments}
This work was partly carried out at the Fields Institute during the Focus Program on Geometric Constraint Systems in 2023. I am grateful to Brigitte Severius for help in translating a tricky passage in Grassmann's papers.

\bibliography{bibl.bib}
\bibliographystyle{plain}

\end{document}